\documentclass[11pt]{amsart}

\usepackage{fullpage}

\numberwithin{equation}{section}

\usepackage{amsmath,amssymb,amscd,amsthm,cite}
\usepackage[mathscr]{eucal}

\newtheorem{theorem}{Theorem}[section]
\newtheorem{proposition}[theorem]{Proposition}

\newtheorem{lemma}[theorem]{Lemma}

\theoremstyle{definition}

\newcommand{\Z}{\mathbb{Z}}
\newcommand{\N}{\mathbb{N}}
\newcommand{\Q}{\mathbb{Q}}
\newcommand{\p}{\partial}

\newcommand{\CP}{\mathbb{CP}}

\newcommand{\<}{\langle}
\renewcommand{\>}{\rangle}
\newcommand{\eps}{\varepsilon}

\newcommand{\DDD}{\mathcal{D}}
\newcommand{\Mbar}{\overline{\mathcal{M}}}

\newcommand{\virt}{{\textup{vir}}}
\DeclareMathOperator{\ev}{ev}
\newcommand{\LL}{\mathbb{L}}
\newcommand{\EE}{\mathbb{E}}

\DeclareMathOperator{\Res}{Res}

\newcommand{\CF}{\mathcal{F}}

\newcommand{\F}[1]{\mathsf{F}^{[#1]}}
\newcommand{\f}[1]{\mathsf{f}^{[#1]}}
\newcommand{\g}{\mathsf{g}}

\renewcommand{\L}[2]{L_{#1,#2}}
\DeclareMathOperator{\Part}{\mathcal{P}}
\DeclareMathOperator{\len}{\ell}
\begin{document}

\title{Multipoint series of Gromov-Witten invariants of $\CP^1$}

\author{E. Getzler, A. Okounkov,  and R. Pandharipande}

\address{Department of Mathematics, Northwestern University, Evanston,
Illinois}

\address{Department of Mathematics, Princeton University,
Princeton, New Jersey}

\address{Department of Mathematics, Princeton University,
Princeton, New Jersey}

\maketitle

\section{Introduction}

The Gromov-Witten theory of $\CP^1$ has an elegant description in
terms of an integrable system called the Toda hierarchy (Eguchi and
Yang \cite{EY}, Eguchi, Hori and Xiong \cite{EHY}, Pandharipande
\cite{P}, Getzler \cite{G} and Okounkov and Pandharipande \cite{OP}).
In this paper, we derive explicit formulas for the multipoint series
of $\CP^1$ in degree 0 from the Toda hierarchy, using the recursions
of the Toda hierarchy (Getzler \cite{G}). Since the Toda equation
determines the higher degree multipoint series of $\CP^1$ from the
degree 0 series, we obtain inductive formulas for the multipoint
series in all degrees, along with explicit formulas for the Hodge
integrals $\int_{\Mbar_{g,n}}\psi_1^{k_1}\dots\psi_n^{k_n}\lambda_g$
and $\int_{\Mbar_{g,n}}\psi_1^{k_1}\dots\psi_n^{k_n}\lambda_{g-1}$.

Let $\Mbar_{g,n,d}=\Mbar_{g,n}(\CP^1,d)$ be the moduli space of stable
genus $g$, $n$-pointed maps to $\CP^1$ of degree $d$.  Let
$\ev_i:\Mbar_{g,n,d}\to\CP^1$ denote the morphism
$$
\ev_i([f:C\to\CP^1,z_1,\dots,z_n]) = f(z_i)
$$
defined by evaluating a stable map $f:C\to\CP^1$ at the $i$th marked
point $z_i\in C$. Let $\psi_i\in H^2(\Mbar_{g,n,d},\Z)$ be the Chern
class of the line bundle over $\Mbar_{g,n,d}$ whose fibre at the point
$$[f:C\to\CP^1,z_1,\dots,z_n]$$ in the moduli space is the cotangent
line $T^*_{z_i}C$.  Let
$$
[\Mbar_{g,n,d}]^\virt \in H_{2(2g-2+2d+n)}(\Mbar_{g,n,d},\Q)
$$
denote the virtual fundamental class.

Let $H\in H^2(\CP^1,\Z)$ be the cohomology class Poincar\'e dual
to the homology  class of a point. Given $k_i,\ell_i\in\N$, define
$$
\<\tau_{k_1,Q}\dots\tau_{k_m,Q}\tau_{\ell_1,P}\dots\tau_{\ell_n,P}\>_{g,d}
 = \int_{[\Mbar_{g,m+n,d}]^\virt} 
\prod_{i=1}^m \psi_i^{k_i} \ev_i^*(H) \cdot
\prod_{i=1}^{n}  \psi_{m+i}^{\ell_i} \in \Q .
$$

The large phase space is the formal affine space with coordinates
$\{s_k,t_k\}_{k\ge0}$. The genus~$g$ Gromov-Witten potential $\CF_g$
of $\CP^1$ is the generating function on the large phase space given
by the formula
$$
\CF_g = \sum_{d=0}^\infty q^d \sum_{m,n=0}^\infty \frac{1}{m!\,n!}
\sum_{k_1,\dots,k_m} \sum_{l_1,\dots,l_n} t_{k_1}\dots t_{k_m}
s_{l_1}\dots s_{l_n}
\<\tau_{k_1,Q}\dots\tau_{k_m,Q}\tau_{l_1,P}\dots\tau_{l_n,P}\>_{g,d} .
$$
The total Gromov-Witten potential is obtained by combining these
potentials into a power series
$$
\CF = \sum_{g=0}^\infty \eps^{2g} \CF_g .
$$
This differs by a factor of $\eps^2$ from the total
Gromov-Witten potential of the physics literature.

Denote the constant vector fields $\p/\p s_k$ and $\p/\p t_k$ on the
large phase space by $\p_{k,P}$ and $\p_{k,Q}$. We will use the
abbreviations $\p$ and $\p_Q$ for $\p_{0,P}$ and $\p_{0,Q}$. We use
the following notation for the partial derivatives of the total
potential:
$$
\<\<\tau_{k_1,Q}\dots\tau_{k_m,Q}\tau_{l_1,P}\dots\tau_{l_n,P}\>\>
= \p_{k_1,Q}\dots\p_{k_m,Q}\p_{l_1,P}\dots\p_{l_n,P}\CF .
$$

We now recall the recursions which characterize the Toda
hierarchy. Define the functions $u=\nabla^2\CF$ and $v=\nabla\p_Q\CF$
on the large phase space, where
\begin{equation}
\nabla  = \eps^{-1}(e^{\eps \partial/2} - e^{-\eps \partial /2}) .
\end{equation}
The total Gromov-Witten potential $\CF$ satisfies the Toda equation
(Okounkov and Pandharipande \cite{OP2})
\begin{equation} \label{toeq}
\p_Q^2\CF = qe^u .
\end{equation}
However, this equation yields no information about the degree 0
potential $\lim_{q\to0}\CF$.

The Toda hierarchy is characterized by the following recursions
(Getzler \cite{G}): if $k$ is a positive integer,
\begin{align} \label{recurse}
\DDD\<\< \tau_{k-1,Q} \>\> &= (k+1) \nabla \<\< \tau_{k,Q} \>\> ,
& \DDD \<\< \tau_{k-1,P} \>\> &= k \nabla \<\< \tau_{k,P} \>\> + 2
\nabla \<\< \tau_{k-1, Q} \>\> ,
\end{align}
where $\DDD$ denotes the operator:
$\DDD=v\nabla+(e^{\eps\p/2}+e^{-\eps\p/2})\p_Q$. The recursions
\eqref{recurse} are proved by combining the partial
Toda hierarchy on the subspace of the large phase space where $s_l=0$
for $l$ (Okounkov and Pandharipande \cite{OP}) with the Virasoro
constraints (Givental \cite{Giv}, Okounkov and Pandharipande
\cite{OP3}).

The multipoint series is closely related to the Gromov-Witten
potential $\CF$ of $\CP^1$. If $d>0$, the degree $d$ multipoint series
$\CF^d[y|z]$ is defined by
$$
\CF^d[y_1,\dots,y_m|z_1,\dots,z_n] =
\sum_{g=0}^\infty \eps^{2g}
\sum_{k_i,l_j} y_1^{k_1} \dots y_m^{k_m} z_1^{l_1} \dots z_n^{l_n} 
\<\tau_{k_1,Q}\dots\tau_{k_m,Q}\tau_{l_1,P}\dots\tau_{l_n,P}\>_{g,d} .
$$
If $m=0$, so that there are no $y$ variables, we write
$\CF^d[z_1,\dots,z_n]$ instead of $\CF^d[~|z_1,\dots,z_n]$.

The multipoint series may also be defined when $d=0$, except that we
modify the definition in two exceptional cases, in order to adjust for
the fact that the moduli spaces $\Mbar_{0,1}$ and $\Mbar_{0,2}$ are empty:
\begin{align*}
\CF^0[y|~] &= y^{-2} + \sum_{g=1}^\infty \eps^{2g}
\sum_k y^k \<\tau_{k,Q}\>_{g,0} , \\ 
\CF^0[y|z] &= (y+z)^{-1} + \sum_{g=1}^\infty \eps^{2g}
\sum_{k,l} y^k z^l \<\tau_{k,Q}\tau_{l,P}\>_{g,0} .
\end{align*}
When $d=0$, the only nonvanishing multipoint series are
$\CF^0[y|z_1,\dots,z_n]$ and $\CF^0[z_1,\dots,z_n]$.

Define auxiliary functions $L_{a,i}(z)$ for $a$ and $i$ in $\N$, by
the formulas
\begin{align*}
\L{a}{0}(z) &= \int_0^z y^{a-1} \biggl( 1 - \frac{\eps y/2}{\tanh(\eps y/2)}
\biggr) dy = - \sum_{k=2}^\infty \frac{\eps^k B_kz^{k+a}}{(k+a)k!} \, , \\
\L{a}{i}(z) &= i \int_0^z \L{a}{i-1}(y) \, dy = - \sum_{k=2}^\infty
\frac{i!\,\eps^kB_kz^{k+a+i}}{(k+a)\dots(k+a+i)k!} \, ,
\end{align*}
where $B_k$ is the $k$th Bernoulli number. 

The following formulas are proved using the recursions \eqref{recurse}.
\begin{theorem} \label{ooo}
Let $Z=z_1+\dots+z_n$, and for $I\subset \{ 1,\dots,n \}$, let
$Z_I=\sum_{i\in I}z_i$. Then
\begin{multline*}
\CF^0[y|z_1,\dots,z_n] =
\biggl(\frac{\eps(y+Z)/2}{\sinh(\eps (y+Z)/2)}\biggr) (y+Z)^{n-2} , \\[15pt]
\shoveleft{\CF^0[z_1,\dots,z_n] = - 2 \biggl(\frac{\eps
Z/2}{\sinh(\eps Z/2)} \biggr) \times} \\
\times Z^{-1} \Biggl( \sum_{\substack{J\subset\{1,\dots,n\}\\|J|<n}}
\sum_{i=1}^{n-|J|} \binom{n-|J|-1}{i-1} (Z-Z_J)^{n-|J|-i}
L_{|J|-1,i}(Z_J) + L_{n-1,0}(Z) \Biggr) .
\end{multline*}
\end{theorem}

There is an overlap between this paper and \cite{GP}, where the
consequences of the Virasoro constraints for $\CP^1$ in degree 0 were
analyzed: we showed that the Virasoro constraints determine the two
degree 0 multipoint series up to overall constants in each genus. The
formula which we prove for $\CF^0[y_1| z_1,\dots,z_n]$ was conjectured
in \cite{GP} and proved, using localization, in \cite{FP1},
\cite{FP2}. Our main new result is the explicit identification of the
more subtle series $\CF^0[ z_1, \dots,z_n]$.

\subsection*{Acknowledgements}
E.~G. and R.~P. thank the organizers of the Workshop on Duality (ITP Santa
Barbara, Summer 2001) for providing the opportunity for collaborating on
this paper. The authors thank T. Eguchi and A. Givental for discussions
about the Toda conjecture and Gromov-Witten theory of $\CP^1$, and
M. Bhargava and K. Kedlaya for help with the proof of Lemma \ref{Lqpp}.

E.~G.\ was partially supported by DMS-0072508, A.~O.\ by DMS-0096246
and fellowships from the Sloan and Packard foundations, and R.~P.\ by
DMS-0071473 and fellowships from the Sloan and Packard foundations.

\section{Hodge series}

Recall the explicit formula for the degree 0 Gromov-Witten
invariants of $\CP^1$ in terms of Hodge integrals \cite{GP}. Let
$\Mbar_{g,n}$ be the moduli space of stable genus $g$, $n$-pointed
curves. Let $$\psi_i=c_1(\LL_i), \ \ \lambda_j = c_j(\EE),$$ where
$\LL_i$ is the $i$th cotangent line bundle and $\EE$ is the Hodge
bundle on $\Mbar_{g,n}$. For $h\ge0$, let
\begin{equation} \label{F}
\F{h} = \sum_{g=h}^\infty (-\eps^2)^g \sum_{n=0}^\infty
\frac{1}{n!} \sum_{k_1,\dots,k_n} s_{k_1}\dots s_{k_n}
\int_{\Mbar_{g,n}} \psi_1^{k_1}\dots\psi_n^{k_n}\lambda_{g-h} .
\end{equation}
\begin{proposition} \label{zero}
The degree 0 Gromov-Witten potential of $\CP^1$ is given by the
formula:
$$
\lim_{q\to0} \CF = \sum_{k=0}^\infty t_k \frac{\p\F{0}}{\p s_k} - 2 \F{1} .
$$
\end{proposition}

Let $f$ be a power series in the variables $s_*=\{s_k\mid k\ge0\}$.
Denote the vector field $\p/\p s_k$ by $\p_k$, and the vector field
$\p_0$ by $\p$. Let $\p(z)$ be the generating function of vector
fields
$$
\p(z) = \sum_{k=0}^\infty z^k \p_k .
$$
If $n>0$, let $f[z_1,\dots,z_n]$ be the generating function
$$
f[z_1,\dots,z_n] = \p(z_1)\dots\p(z_n) f|_{s_*=0} .
$$
Proposition \ref{zero} yields the following expression for the
generating functions of degree 0 Gromov-Witten invariants of $\CP^1$:
\begin{align} \label{fff}
\CF^0[y_1|z_1,\dots,z_n] &= \F{0}[y_1,z_1,\dots,z_n] , &
\CF^{0}[z_1,\dots,z_n] &= -2 \F{1}[z_1,\dots,z_n] .
\end{align}

The string equation in Gromov-Witten theory implies the following
differential equation for the Hodge series:
\begin{equation} \label{string}
\p\F{h} = \delta_{h,0} \biggl( \frac{s_0^2}{2} - \frac{\eps^2}{24}
\biggr) + \sum_{n=0}^\infty s_{n+1} \p_n\F{h} .
\end{equation}
We obtain the following relations between the series
$\F{h}[z_1,\dots,z_n]$:
\begin{equation} \label{String}
\F{h}[z_1,\dots,z_n,0] = Z \F{h}[z_1,\dots,z_n] + \delta_{h,0}
\delta_{n,2} .
\end{equation}

Let $\f{h}$ be the generating function $\f{h}=\nabla\F{h}$. The
string equation \eqref{string} for $\F{h}$ implies the string
equation for $\f{h}$:
\begin{equation} \label{fstring}
\p\f{h} = \delta_{h,0} s_0 + \sum_{n=0}^\infty s_{n+1} \p_n\f{h} .
\end{equation}
Although we are really more interested in the series
$\F{h}[z_1,\dots,z_n]$, the recursions \eqref{recurse} lead more
naturally to formulas for the series $\f{h}[z_1,\dots,z_n]$.
The string equation gives a formula relating these two
sets of series:
\begin{proposition} \label{fF}
If $n\ge3$, then
$$
\F{h}[z_1,\dots,z_n] = \biggl(\frac{\eps Z/2}{\sinh(\eps Z/2)} \biggr)
Z^{-1} \f{h}[z_1,\dots,z_n] .
$$
\end{proposition}
\begin{proof}
We have
$$
\f{h}[z_1,\dots,z_n] = \sum_{g=h}^\infty \sum_{m=0}^\infty
\frac{(-\eps^2)^g\eps^{2m}}{2^{2m}(2m+1)!} \int_{\Mbar_{g,n+2m+1}}
\frac{\lambda_{g-h}}{(1-z_1\psi_1)\dots(1-z_n\psi_n)} .
$$
By the string equation \eqref{string},
$$
\int_{\Mbar_{g,n+2m+1}}
\frac{\lambda_{g-h}}{(1-z_1\psi_1)\dots(1-z_n\psi_n)}
= Z^{2m+1} \int_{\Mbar_{g,n}}
\frac{\lambda_{g-h}}{(1-z_1\psi_1)\dots(1-z_n\psi_n)} ,
$$
and the result follows on summing over $m$.
\end{proof}

\section{Calculation of $\f{0}$}

We will now derive a formula for the Hodge series
$\f{0}[z_1,\dots,z_n]$. We actually establish a more general result,
for the generating function $e^{x\p}\f{0}$, obtained from $\f{0}$ by
translating $s_0$ by $x$.
\begin{theorem} \label{f0}
If $n>1$, $(e^{x\p}\f{0})[z_1,\dots,z_n]=Z^{n-2}e^{xZ}$.
\end{theorem}
\begin{proof}
By Proposition \ref{zero}, we have
$$
\text{$\displaystyle\lim_{q\to0} \nabla \<\<\tau_{n,Q}\>\>=\p_n\f{0}$
and $\displaystyle\lim_{q\to0} \p_Q \<\<\tau_{n,Q}\>\> = 0$.}
$$
In particular, $\lim_{q\to0}v=\p\f{0}$. Since we work here
exclusively in degree 0, we will denote $\p\f{0}$ by $v$. The
$q\to 0$ limit of the recursion
$$\DDD\<\< \tau_{n-1,Q} \>\> = (n+1) \nabla
\<\< \tau_{n,Q}\>\> $$ yields the following equation:
$v \p_{n-1}\f{0} = (n+1) \p_n\f{0}$. Thus, we find that
\begin{equation} \label{pnv}
\p_n\f{0} = \frac{v^{n+1}}{(n+1)!} .
\end{equation}
Summing over $n$, we see that
\begin{equation} \label{pf}
\p(z)\f{0} = z^{-1}(e^{zv}-1) .
\end{equation}

We now prove by induction that for $n>1$,
$(n-1)!\,\p(z_1)\dots\p(z_n)\f{0}$ equals the residue at $x=0$ of
$$
x^{-n} Z^{-1} \exp(Ze^{x\p}v) .
$$
Replacing $z$ by $z_1$ in \eqref{pf} and applying the operator
$\p(z_2)$ shows that
$$
\p(z_1)\p(z_2)\f{0} = (\p(z_2)v) \exp(z_1v) = \p v \exp((z_1+z_2)v) ,
$$
giving the case $n=2$. The induction step is as follows:
\begin{align*}
(n-1)! \, \p(z_1)\dots\p(z_{n+1})\f{0} &= \p(z_{n+1}) \Res_0( x^{-n}
Z^{-1} \exp(Ze^{x\p}v) ) \\
&= \Res_0( x^{-n} (\p(z_{n+1})e^{x\p}v) \exp(Ze^{x\p}v) ) \\
&= \Res_0( x^{-n} (e^{x\p}\p v) \exp((Z+z_{n+1})e^{x\p}v) ) \\
&= \Res_0 \biggl( x^{-n} \frac{\p}{\p x}
\frac{\exp((Z+z_{n+1})e^{x\p}v)}{Z+z_{n+1}} \biggr) .
\end{align*}

By the string equation \eqref{fstring}, $(e^{x\p}v)|_{s_*=0}=x$.
We then conclude for $n>1$,
\begin{align*}
(n-1)!\,(e^{y\p}\f{0})[z_1,\dots,z_n] &= \Res_{x=0}( x^{-n} Z^{-1}
\exp((x+y)Z) ) \\ &= (n-1)! \, Z^{n-2} e^{yZ} .
\qedhere\end{align*}
\end{proof}

The formula for $\CF^0[y|z_1,\dots,z_n]$ in Theorem \ref{ooo} is an
immediate consequence of \eqref{lplp}, \eqref{fff} and \eqref{String}
(which is needed to prove the cases $n=0$ and $n=1$).

Theorem \ref{f0} is closely related to results of the paper
\cite{FP2}. By Proposition~\ref{fF}, we see that for $n\ge3$,
\begin{equation}
\label{lplp}
\F{0}[z_1, \dots, z_n] = \biggl(\frac{\eps Z/2}{\sinh(\eps Z/2)}
\biggr) Z^{n-3} .
\end{equation}
Using \eqref{String}, we may also handle the cases with $n<3$. For
example, since
$$
\F{0}[z] = z^{-2} \biggl(\frac{\eps z/2}{\sinh(\eps z/2)} - 1 \biggr) ,
$$
we recover a formula proved in \cite{FP1}: for $g>0$,
$$
\int_{\Mbar_{g,1}} \psi_1^{2g-2} \lambda_g = (-1)^g \bigl( 2^{1-2g} - 1 \bigr)
\frac{B_{2g}}{(2g)!} .
$$
We do not know how to derive the above formula directly from the Virasoro
constraints, even though in principle, the Virasoro conjecture determines
$\CF$ from the genus $0$ Gromov-Witten potential $\CF_0$ (Dubrovin and
Zhang \cite{DZ}).

\section{Calculation of $\f{1}$}

In this section, we prove the formula of Theorem \ref{ooo} for
$\CF[z_1,\dots,z_n]$. We actually work with the series
$\f{1}[z_1,\dots,z_n]$, establishing a recursion which determines
$\f{1}[z_1,\dots,z_n]$ from $\f{1}[z_1,\dots,z_m]$, $m<n$. The formula
for $\CF[z_1,\dots,z_n]$ follows on application of Proposition
\ref{fF} and \eqref{fff}.
\begin{theorem} \label{bbl}
$$
d\f{1}[z_1,\dots,z_n] = \sum_{\substack{I\subset\{1,\dots,n\}\\|I|<n}}
(Z-Z_I)^{n-|I|-1} \f{1}[z_I] \, dZ_I + Z^{n-1} d \log\biggl(
\frac{\eps Z/2}{\sinh(\eps Z/2)} \biggr)
$$
\end{theorem}
\begin{proof}
By \eqref{recurse}, we see 
\begin{equation} \label{Recurse}
\sum_{n=0}^\infty z^n \DDD\<\<\tau_{n,P}\>\> = \frac{\p}{\p z}
\sum_{n=0}^\infty z^n \nabla\<\<\tau_{n,P}\>\> + 2
\sum_{n=0}^\infty z^n \nabla\<\<\tau_{n,Q}\>\> .
\end{equation}
From Proposition \ref{zero}, we obtain the formulas:
\begin{align*}
\lim_{q\to0} \sum_{n=0}^\infty z^n \nabla\<\<\tau_{n,P}\>\>
&= \sum_{k=0}^\infty t_k\p_k\p(z)\f{0} - 2\p(z)\f{1}
= e^{zv} \sum_{k=0}^\infty t_k\p_kv - 2\p(z)\f{1} , \\
\lim_{q\to0} \sum_{n=0}^\infty z^n \nabla\<\<\tau_{n,Q}\>\> &=
\p(z)\f{0} , \\
\lim_{q\to0} \sum_{n=0}^\infty z^n \p_Q\<\<\tau_{n,P}\>\> &= \p(z)
\biggl( \frac{\eps\p/2}{\sinh(\eps\p/2)} \biggr) \f{0} .
\end{align*}
We therefore find:
\begin{align*}
\lim_{q\to0} \sum_{n=0}^\infty z^n \DDD\<\<\tau_{n,P}\>\> &=
\lim_{q\to0} \sum_{n=0}^\infty z^n \bigl( v\nabla\<\<\tau_{n,P}\>\> +
(e^{\eps\p/2} + e^{-\eps\p/2}) \p_Q \<\<\tau_{n,P}\>\> \bigr) \\
&= v e^{zv} \sum_{k=0}^\infty t_k\p_kv - 2v\p(z)\f{1} + \p(z)
\biggl( \frac{\eps\p}{\tanh(\eps\p/2)} \biggr) \f{0} .
\end{align*}
The terms involving the variables $t_k$ cancel in the recursion
\eqref{Recurse}, and we are left with the formula:
\begin{align*}
\frac{\p}{\p z} (\p(z)\f{1}) &= v\p(z)\f{1} + \biggl( 1 -
\frac{\p/2}{\tanh(\p/2)} \biggr) \p(z) \f{0} \\
&= v\p(z)\f{1} + \biggl( 1 - \frac{\eps\p_x/2}{\tanh(\eps\p_x/2)}
\biggr)_{x=0} \p(z) e^{x\p}\f{0} .
\end{align*}
Substituting $z_i$ for $z$ in this formula and applying the
operator $\p(z_1)\dots\widehat{\p(z_i)}\dots\p(z_n)$, we obtain:
\begin{multline*}
\frac{\p}{\p z_i} \bigl( \p(z_1)\dots\p(z_n)\f{1} \bigr) \\
= \sum_{\substack{I\subset\{1,\dots,n\}\\i\in I}}
\prod_{j\notin I} \p(z_j)v \, \prod_{j\in I} \p(z_j)\f{1} +
\biggl( 1 - \frac{\eps\p_x/2}{\tanh(\eps\p_x/2)} \biggr)_{x=0}
\p(z_1)\dots\p(z_n) e^{x\p} \f{0} .
\end{multline*}
Note that for $n>0$,
$v[z_1,\dots,z_n]=\p\f{0}[z_1,\dots,z_n]=Z^{n-1}$. Hence, taking
$s_*=0$, we see:
\begin{align*}
\frac{\p}{\p z_i} \f{1}[z_1,\dots,z_n]
&= \sum_{\substack{I\subset\{1,\dots,n\}\\i\in I}} (Z-Z_I)^{n-|I|-1}
\f{1}[z_I] + \biggl( 1 - \frac{\eps\p_x/2}{\tanh(\eps\p_x/2)}
\biggr)_{x=0} Z^{n-2} e^{xZ} \\
&= \sum_{\substack{I\subset\{1,\dots,n\}\\i\in I}} (Z-Z_I)^{n-|I|-1}
\f{1}[z_I] + Z^{n-2} \biggl( 1 - \frac{\eps Z/2}{\tanh(\eps Z/2)}
\biggr) .
\end{align*}
Multiplying by $dz_i$, summing over $i$, and using the identity
$$
Z^{n-2} \biggl( 1 - \frac{\eps Z/2}{\tanh(\eps Z/2)} \biggr) dZ =
Z^{n-1} d \log\biggl( \frac{\eps Z/2}{\sinh(\eps Z/2)} \biggr) ,
$$
the theorem is proved.
\end{proof}

The string equation \eqref{fstring} implies $\f{1}[0,\ldots,0]=0$.
Hence, Theorem \ref{bbl} determines the functions
$\f{1}[z_1,\ldots,z_n]$. For example,
$$
\f{1}[z] = \log\biggl( \frac{\eps z/2}{\sinh(\eps z/2)} \biggr) = -
\sum_{k=2}^\infty \frac{\eps^kB_k}{k!} \frac{z^k}{k},
$$
which is equivalent to the formula for $\F{1}[z]$ in \cite{FP1},
by Proposition \ref{fF}.

Starting with the above formula for $\f{1}[z]$, Theorem \ref{bbl}
yields:
$$
\f{1}[z_1,z_2] = - \sum_{k=2}^\infty \frac{\eps^kB_k}{k!} \biggl(
\frac{(z_1+z_2)^{k+1}}{k+1} + \frac{z_1^{k+1}+z_2^{k+1}}{k(k+1)}
\biggr) .
$$
We now find a general formula for $\f{1}[z_1,\dots,z_n]$.
Let
$$
\g[z_1,\dots,z_n] = \sum_{\substack{J\subset\{1,\dots,n\}\\|J|<n}}
\sum_{i=1}^{n-|J|} \binom{n-|J|-1}{i-1} (Z-Z_J)^{n-|J|-i}
L_{|J|-1,i}(Z_J) + L_{n-1,0}(Z) .
$$
\begin{theorem}
$\f{1}[z_1,\dots,z_n]=\g[z_1,\dots,z_n]$.
\end{theorem}
\begin{proof}
For $n=1$, both $\f{1}[z_1]$ and $\g[z_1]$ 
equal $L_{0,0}(Z)$. We prove the Theorem by
showing the functions $\g[z_1,\dots,z_n]$ satisfy the
recursion of Theorem \ref{bbl}, which may be restated as:
\begin{equation} \label{r}
d\g[z_1,\dots,z_n] - dL_{n-1,0}(Z) =
\sum_{\substack{I\subset\{1,\dots,n\}\\|I|<n}} (Z-Z_I)^{n-|I|-1}
\g[z_I] \, dZ_I ,
\end{equation}
where $\g[z_I]= \g[z_{i_1}, \dots, z_{i_|I|}]$ for $I= \{i_1< \dots 
< i_{|I|} \}$.
The left-hand side of \eqref{r} equals:
\begin{equation} \label{s}
\sum_{\substack{J\subset\{1,\dots,n\}\\|J|<n}}
\sum_{i=0}^{n-|J|} \binom{n-|J|-1}{i} (Z-Z_J)^{n-|J|-i-1}
L_{|J|-1,i}(Z_J) ( (i+1) dZ_J + i d(Z-Z_J) ) .
\end{equation}

Lemma \ref{Lqpp} below will be used to handle \eqref{s}; the first of
the identities is due to Hurwitz \cite{hur}, while the second is
closely related. (See also Pitman \cite{for}.)
\begin{lemma} \label{Lqpp}
\begin{align*}
(x+y)(x+y+Z)^{n-1} &= xy \sum_{I\subset\{1,\dots,n\}} (x+Z_I)^{|I|-1}
(y+Z-Z_I)^{n-|I|-1} \\
(x+y+Z)^{n-1} \, dZ &= y \sum_{I\subset\{1,\dots,n\}} (x+Z_I)^{|I|-1}
(y+Z-Z_I)^{n-|I|-1} \, dZ_I
\end{align*}
\end{lemma}
\begin{proof}
We require a result from the theory of tree enumeration. A tree $T$ with
vertices $\{1,\dots,n\}$ is a subset of the set of pairs $\{(i,j)\mid 1\le
i<j\le n\}$ of cardinality $n-1$ such that the associated graph is
connected. Let $d(T,i)$ be the valence of the vertex $i$ in the tree $T$. A
basic combinatoric identity is:
$$
Z^{n-2} = \sum_{\substack{\text{trees $T$ with vertices} \\
\{1,\dots,n\}}} \prod_{i=1}^n z_i^{d(T,i)-1} .
$$

Thus, the summand of the right-hand side of the first identity is obtained
by summing over pairs of trees, one of which has vertices $\{x\}\cup I$,
and the other of which has vertices $\{y\}\cup\{1,\dots,n\}\backslash
I$. Similarly, the coefficient of $dz_k$ in the summand of the second
identity is obtained by summing over pairs of trees, one of which has
vertices $\{x\}\cup I$, and the other of which has vertices
$\{y\}\cup\{1,\dots,n\}\backslash I$, subject to the additional condition
that $k\in I$.

Let us prove the first formula. The pairs of trees contributing to the
right-hand side may be obtained by taking a tree with vertices
$\{w\}\cup\{1,\dots,n\}$, and cutting the vertex $w$ into two vertices $x$
and $y$ of valence $e$ and $d-e$. If the vertex $w$ had valence $d$, and
hence weight $w^{d-1}$, then after being split, the contribution
$x^{e-1}y^{d-e-1}$ is obtained.  Summation over $e$ yields $(xy)^{-1}
(x+y)^d$. The left-hand side of the identity is the result of performing
the replacement of $w^{d-1}$ by $(xy)^{-1} (x+y)^d$ on the sum over trees
$(w+Z)^{n-1}$.

The proof of the second formula is similar, except that now, we insist that
the $e$ edges emerging from $x$ include the distinguished edge emerging
from $w$ which lies on the unique path from $w$ to $k$. Now, summing over
$e$ yields $y^{-1} (x+y)^{d-1}$. The left-hand side of the second identity
is the result of performing the replacement of $w^{d-1}$ by $y^{-1}
(x+y)^{d-1}$ on the sum over trees $(w+Z)^{n-1}$.
\end{proof}

Taking a derivative of both sides of the first identity of Lemma \ref{Lqpp}
with respect to $y$ and setting $y=0$ yields
\begin{align*}
(x+Z)^{n-1} + (n-1) x(x+Z)^{n-2} &= x
\sum_{\substack{I\subset\{1,\dots,n\}\\|I|<n}} (x+Z_I)^{|I|-1}
(Z-Z_I)^{n-|I|-1} \\
&= x \sum_{\substack{I\subset\{1,\dots,n\}\\0<|I|<n}} (x+Z_I)^{|I|-1}
(Z-Z_I)^{n-|I|-1} + Z^{n-1}
\end{align*}
Extracting the coefficient of $x^i$, $i>0$, we see that
\begin{equation} \label{lqpp}
(i+1) \binom{n-1}{i} Z^{n-i-1} =
\sum_{\substack{I\subset\{1,\dots,n\}\\0<|I|<n}}
\binom{|I|-1}{i-1} Z_I^{|I|-i} (Z-Z_I)^{n-|I|-1} .
\end{equation}
Likewise, taking a derivative of both sides of the second identity
with respect to $y$ and setting $y=0$ yields
$$
(n-1) (x+Z)^{n-2} \, dZ = \sum_{\substack{I\subset\{1,\dots,n\}\\
|I|<n}} (x+Z_I)^{|I|-1} (Z-Z_I)^{n-|I|-1} \, dZ_I .
$$
Extracting the coefficient of $x^{i-1}$ gives
\begin{equation} \label{lqppp}
i \binom{n-1}{i} Z^{n-i-1} \, dZ =
\sum_{\substack{I\subset\{1,\dots,n\}\\0<|I|<n}}
\binom{|I|-1}{i-1} Z_I^{|I|-i} (Z-Z_I)^{n-|I|-1} \, dZ_I .
\end{equation}
Using formulas \eqref{lqpp} and \eqref{lqppp}, we easily see that
\eqref{s} equals the right-hand side of \eqref{r}.
\end{proof}

\section{Higher degree multipoint series}

In this section, we derive an inductive formula from the Toda equation
\eqref{toeq} for the higher degree multipoint series. In particular,
this will show that the multipoint series of the Gromov-Witten theory
of $\CP^1$ are integrals of trigonometric functions in all degrees $d$.

From the Toda equation \eqref{toeq}, we see that if $d>0$,
$$
\p_Q^2\CF^d = \sum_{k=1}^\infty \frac{1}{k!} \sum_{d_1+\dots+d_k=d-1}
\nabla^2\CF^{d_1} \dots \nabla^2\CF^{d_k} .
$$
Using the divisor equation
$$
\p_Q\CF^d = d \cdot \CF^d + \sum_{k=0}^\infty s_{k+1} \p_{k,Q}\CF^d ,
$$
we see that
$$
\p_Q^2\CF^d = d^2 \cdot \CF + 2d \sum_{k=0}^\infty s_{k+1} \p_{k,Q}\CF^d +
\sum_{j,k=0}^\infty s_{j+1} s_{k+1} \p_{j,Q}\p_{k,Q}\CF^d .
$$
We conclude that
\begin{equation} \label{toda}
\CF = - \frac{2}{d} \sum_{k=0}^\infty s_{k+1} \p_{k,Q}\CF^d - \frac{1}{d^2}
\sum_{j,k=0}^\infty s_{j+1} s_{k+1} \p_{j,Q}\p_{k,Q}\CF^d
+ \sum_{\ell=1}^\infty \frac{1}{\ell!} \sum_{d_1+\dots+d_\ell=d-1}
\nabla^2\CF^{d_1} \dots \nabla^2\CF^{d_\ell} .
\end{equation}

Let $\Part_{n,m}$ be the set of partitions with nonempty parts of the set
$\{y_1,\dots,y_m,z_1,\dots,z_n\}$. For $P\in\Part_{n,m}$, let
$\{P_1,\dots,P_{\len(P)}\}$ denote the parts of $P$. Let $y_{P_i}$ and
$z_{P_i}$ denote the $y$ and $z$ variables of the part $P_i$. Let $Y_{P_i}$
and $Z_{P_i}$ denote the sums of the variable sets $y_{P_i}$ and $z_{P_i}$
respectively. We see from \eqref{toda} that for $d>0$,
\begin{multline} \label{recct}
\CF^d[y_1,\dots,y_m|z_1,\dots,z_n] = \\
\begin{aligned}
{}& -\frac{2}{d} \sum_{1\leq i\leq n} z_i
\CF^d[y_1,\dots,y_m,z_i|z_1,\dots,\hat{z_i},\dots z_n] \\
{}& -\frac{2}{d^2}\sum_{1\leq  i<j \leq n} z_i z_j
\CF^d[y_1,\dots,y_m,z_i,z_j|z_1,\dots,\hat{z_i},\dots,\hat{z_j},\dots z_m] \\
{}& + \frac{1}{d^2} \sum_{P\in\Part_{m,n}}
\sum_{d_1+\dots+d_{\len(P)}=d-1} \prod_{i=1}^{\len(P)}
\biggl( \frac{\sinh(\eps(Y_{P_i}+Z_{P_i})/2)}{\eps(Y_{P_i}+Z_{P_i})/2}
\biggr)^2 \CF^{d_i}[y_{P_i}|z_{P_i},0,0] .
\end{aligned}
\end{multline}
This equation, together with the formulas of Theorem \ref{ooo} for the
degree 0 multipoint series, gives explicit formulas for the higher degree
multipoint series.

\end{document}